\documentstyle{amsppt}
\magnification=\magstep1
\loadbold
\baselineskip=18pt

\topmatter
\title
{On the Dirichlet problem of Landau-Lifshitz-Maxwell equations}
\endtitle
\leftheadtext{ Landau-Lifshitz-Maxwell equations }
\author
Jian Zhai
\endauthor
\affil Department of Mathematics, Zhejiang University, Hangzhou
310027, P.R.China; Email: jzhai\@zju.edu.cn
\endaffil
\rightheadtext{J. Zhai}
\abstract We prove the  existence and uniqueness of non-trivial
stable solutions to Landau-Lifshitz-Maxwell equations with
Dirichlet boundary condition for large anisotropies and small
domains, where the domains are non-simply connected.
\endabstract
\thanks
\endthanks
\endtopmatter

\bigskip

\indent
{\bf 1 Introduction}

We seek for a solution $u=(u_1,u_2,u_3): \Omega \to S^2\subset\Bbb R^3$ of the Landau-Lifshitz-Maxwell equations with the Dirichlet boundary condition
$$\left\{\aligned
&\Delta u+|\nabla u|^2 u +H-(H\cdot u)u-\lambda
(W_u(u)-(W_u(u)\cdot u)u)=0\quad\text{in}\quad \Omega\\
&u=g\quad\text{on}\quad \partial\Omega\\
&\text{curl H}=0,\quad \text{div}(H+u\chi_\Omega)=0,\quad\text{in}\quad \Bbb R^3
\endaligned\right.
\tag1.1
$$
where $\Omega$ is a non-simply connected bounded domain in $\Bbb R^3$ with uniformly $C^4$ boundary, $\lambda > 0$ is a parameter, $g \in C^{3+\alpha_0}(\partial\Omega,
S^2\cap\{u_3=0\})$ $(0<\alpha_0<1)$ and $\chi_\Omega$ is the
characteristic function of the domain $\overline\Omega$
$$
\chi_\Omega(x)=\left\{\aligned
                      &1,\quad\forall x\in\overline\Omega\\
                      &0,\quad\forall x\not\in\overline\Omega.
                       \endaligned\right.
$$
  For thin films with in plane magnetization (c.f. [HS]), as a first approximation, we can assume $W(u)=u_3^2$ and denote $W_u(u)=(0,0,2u_3)$.

$H$ is the demagnetizing field generated by the magnetization $u$
and determined by the Maxwell's equation. (1.1) is the
Euler-Lagrange equation of the Landau-Lifshitz energy functional
$$E_\lambda(u)=\int_\Omega \frac1{2}|\nabla u|^2-\frac12 u\cdot H +\lambda
W(u)dx\tag1.2
$$
on  $H^1_g(\Omega,S^2)$. Here
$$
H^1_g(\Omega,S^2):=\{u-g\in H^1_0(\Omega,\Bbb R^3)|\,\,u(\Omega)\subset S^2\}
$$
is well defined for $\Omega$ and $g$.

Functional (1.2) was first derived for ferromagnetic problem by
Landau and Lifshitz [LL] in 1935.  The equation (1.1) is  the static
equivalent of the time-dependent Landau-Lifshitz-Maxwell equations
(c.f.[Z3]-[Z11])
$$
\left\{\aligned
\frac{\partial u}{\partial t} &= -u\times(u\times(\Delta u-\lambda
W_u(u)+H))\\
& +\gamma u\times (\Delta u-\lambda W_u(u)+H) \quad\text{in}\quad \Omega\times(0,\infty)\\
 u &=g \quad\text{on}\quad
\partial\Omega\times(0,\infty)\\
u &\in S^2 \quad\text{in}\quad \Omega\times(0,\infty)\\
&\text{curl}\, H=0,\quad\text{div} (H+u\chi_\Omega)=0,\quad\text{in}\quad\Bbb R^3.
\endaligned\right.\tag1.3
$$
Here $\gamma\geqq 0$ is called the damping constant and small
$\gamma$ means large damping.

Afterhere, a solution $u(x)$ will  be expressed in
$$u(x)=(\cos\xi(x)\cos\theta(x),\cos\xi(x)\sin\theta(x),\sin\xi(x)),\quad (-\frac\pi2\leqq \xi \leqq \frac\pi2,\,\, \theta \in S^1=\Bbb
R/2\pi\Bbb Z),
$$
and
$$
g(x)=(\cos\theta_g(x),\sin\theta_g(x),0).
$$
For convenience, we also use the notations
$$u(\theta,\xi)=(\cos\xi(x)\cos\theta(x),\cos\xi(x)\sin\theta(x),\sin\xi(x))$$
 and
$$
u_\theta(\theta,\xi)=(-\cos\xi\sin\theta,\cos\xi\cos\theta,0),\quad u_\xi(\theta,\xi)=(-\sin\xi\cos\theta,-\sin\xi\sin\theta,\cos\xi)
$$
and etc.

In [Z3], the existence and non-existence of non-trivial stable
solutions to Landau-Lifshitz equation (the first equation of (1.1)
with constant $H$) with the Dirichlet boundary condition were
obtained. This paper is the continuation of [Z3]. In this paper,
we study the Landau-Lifshitz-Maxwell equations (1.1) in a
non-simply connected bounded domain of $\Bbb R^3$. The existence
and uniqueness of non-trivial solutions and their stability are
obtained. Precisely, we have

\bigskip

\proclaim{Theorem 1}  There exists $d_0>0$ such that for any
non-simply connected domain $\Omega$ if $|\Omega|\leq d_0$, then there is
$\lambda_0>0$ such that for $\lambda>\lambda_0$, there exists a unique
solution $(u_\lambda, H_\lambda)$
$$\aligned
&u_\lambda(x)=(\cos\xi_\lambda(x)\cos\theta_\lambda(x),\cos\xi_\lambda(x)\sin\theta_\lambda(x),\sin\xi_\lambda(x)) \in C^{2+\alpha}(\overline\Omega)\\
&H_\lambda\in C^{1+\alpha}(\overline\Omega)\cap (\text{BMO}\cap L^{p})(\Bbb R^3), \quad \forall p\in (1,\infty)
\endaligned$$
$(0<\alpha<1)$ to (1.1) corresponding to the homotopy class $[\theta_g]$
of continuous maps from $\Omega$ to $S^2\cap\{u_3=0\}$.  Moreover,
$\theta_\lambda$ is homotopic to $\theta_g$ and
$$
\|\xi_\lambda\|_{C^\alpha(\overline\Omega)} \leqq \frac C{\sqrt{\lambda}},
$$
where the constant $C$ is independent of $\lambda$.

\endproclaim

\bigskip

 A steady state solution
$u_\lambda(x)$ is called exponentially asymptotically stable (c.f.[L]) if  there are $\mu_0>0$ and
$\epsilon, C>0$ such that for all $\bar u\in W^{2,2}(\Omega, S^2)$, if $\|\bar
u-u_\lambda\|_{W^{2,2}(\Omega)}\leq \epsilon$,  then there exists a unique global  solution $(u(x,t),H)$  of (1.3) with the initial data
$u(x,0)=\bar u(x)$, and
$$
\|u(x,t)-u_\lambda(x)\|_{W^{2,2}(\Omega)}\leq Ce^{-\mu_0 t}\|\bar u(x)-u_\lambda(x)\|_{W^{2,2}(\Omega)},\quad\forall t\geq0.
$$

\bigskip

\proclaim{Theorem 2}  Assume that
$\Omega$ is not simply connected and $|\Omega|\leq d_0$ as in Theorem 1. Then there exists a $\lambda_0>0$
and
for $\lambda \geqq \lambda_0$, there exists a $\gamma_0
=\gamma_0(\lambda)>0$ such that, for $\lambda\geqq\lambda_0$ and
$\gamma \in [0,\gamma_0]$,  the solution
$(u_\lambda,H_\lambda)$
obtained in Theorem 1
is exponentially asymptotically stable steady state solutions of the time-dependent Landau-Lifshitz-Maxwell equations (1.3).
\endproclaim

\bigskip

 To prove the theorems, we first solve the Maxwell's equation  and  express $H$
in $u$ and a Calder\'{o}n-Zygmund operator. Then we consider the
limit case of $\lambda\to\infty$ and get a solution in the
homotopy class $[\theta_g]$  of continuous mappings from
$\overline\Omega$ to $S^1$ by Schauder fixed point theorem. We
search solutions for large $\lambda$ in the neighborhood of the
limit case solutions by the Schauder fixed point theorem again. In
the last we analyze the spectrum of the linearized operator in a
detailed way by using the Kato's perturbation theory (c.f. [K])
and to prove that the solutions are exponentially asymptotically
stable steady state solutions of (1.3) by using nonlinear
parabolic equations theory (c.f.[L]).

In [Z11], we proved the existence of non-trivial stable  solutions
to Landau-Lifshitz-Maxwell equations with Neumann boundary
condition for large anisotropies and small domains that are non-simply connected and rotationally invariant around
an axis.

Remark that in the Theorem 1 and Theorem 2 of this paper we do not need the assumption that the domain $\Omega$ is
rotation invariant  as in [Z11] for Neumann
boundary condition. Moreover the assumption of diam$(\Omega)\leq d_0$ used in [Z11] is replaced by the volume $|\Omega|\leq  d_0$.

On the other hand, in [Z3] and [Z11], the key estimate
$$
\|\xi\|_{C^{\alpha_0}(\overline\Omega)}\leq \frac C{\lambda}
$$
 was obtained by using the Campanato inequality. But for the Dirichlet problem considered in this paper, the
Campanato inequality can not be used, because the right of (5.6)
may be non-zero on the boundary of $\Omega$ (c.f. [C]).  Here we
shall apply [L] Theorem 3.1.3 to (5.6) to get a similar estimate
$$
\|\xi\|_{C^{\alpha_0}(\overline\Omega)}\leq \frac
C{\sqrt{\lambda}}.
$$

Related works can also be founded in [JK], [V], [ABV], [KS], [GSh], [SC], [BTW], [RS], [FR], [GS], [HL] and [Z1]-[Z11] and the
papers cited in there.   Visintin [V] proved the existence of a
kind of weak solutions to  Landau-Lifshitz-Maxwell equations. In
two dimensional case, the existence, uniqueness and partial
regularity of (1.3) were considered in  [GS]. Also in the
two dimension, Gustafson and Shatah studied the existence and
stability of localized periodic solution to the Landau-Lifshitz
equation, and Chang, Shatah, Uhlenbeck studied the well-posedness
of the Cauchy problem for Schr\"{o}dinger maps ( there are not $H$
and the first term in the right of the first equation of (1.3)).
In three dimension,  Ball, Taheri and Winter constructed  local
energy minimizers around a fixed constant solution to the model of
micromagnetics.  Their results are different from mine.

This paper consists of five sections. The Maxwell equation is studied in section 2. In section 3, the Landau-Lifshitz-Maxwell equations are expressed in spherical coordinates.   In section 4, we consider the limit case  of $\lambda\to\infty$.   The main theorems of this paper  are proved in section 5-6.

\bigskip

\indent
{\bf 2. Maxwell equations}

The magnetization $u$ and the demagnetizing field $H$ are related by the Maxwell's equation
$$\left\{\aligned
          &\text{curl}H=0\quad\text{in}\quad\Bbb R^3\\
          &\text{div}(H+u\chi_\Omega)=0\quad\text{in}\quad\Bbb R^3.
         \endaligned\right.
\tag2.1$$

First we recall a lemma  proved in [Z11]. For reader's convenience, we also give its proof here.

\proclaim{Lemma 2.1}There exists a continuous linear map $\Cal{L}$
$$
L^2(\Omega;\Bbb R^3)\ni u\mapsto \nabla v,\quad v\in V=\{v\in H^1_{loc}(\Bbb R^3;\Bbb R):\quad \nabla v\in L^2(\Bbb R^3;\Bbb R^3),\,\, \int_\Omega vdx=0\}
$$
such that $H=\nabla v=\Cal{L}( u)$ is the unique solution of  (2.1) in $V$. Moreover,
$\Cal{L}$ is bounded from
$$\aligned
&\Cal{H}^1\text{(Hardy)}\quad\text{to}\quad L^1(\Bbb R^3),\\
&L^\infty(\Omega) \quad\text{to}\quad \text{BMO},\\
&L^p(\Omega)\quad\text{to}\quad L^p(\Bbb R^3), \quad\forall p\in(1,\infty).
\endaligned$$

\endproclaim

$Proof.$ From [JK], we know that there exists a continuous linear map
$$
L^2(\Omega;\Bbb R^3)\ni u\mapsto v\in \{v\in H^1_{loc}(\Bbb R^3;\Bbb R):\quad \nabla v\in L^2(\Bbb R^3;\Bbb R^3),\,\, \int_\Omega vdx=0\}
$$
such that $H=\nabla v=\Cal{L}( u)$ is the unique solution of  (2.1) in $V$.

Note that from
$$
\text{div}(\nabla v+\chi_\Omega u)=0,\quad\text{in}\quad \Bbb R^3,
$$
for any domain $D\supset\supset\Omega$,
$$
0=\int_D\text{div}(\nabla v+\chi_\Omega u)dy=\int_{\partial
D}\frac{\partial}{\partial\nu}v(y)dS(y).
$$
Then for any $x\in \Bbb R^3$, for any ball
$B_R(x)\supset\supset\Omega$, by the Green's representation
formula
$$\aligned
v(x)&=\int_{\partial B_R(x)}v(y)\frac{\partial}{\partial \nu}\frac1{4\pi|x-y|}-\frac1{4\pi|x-y|}\frac{\partial}{\partial \nu}v(y)dS(y)\\
&-\int_{B_R}\frac1{4\pi|x-y|}\text{div}(\chi_\Omega(y)u(y))dy\\
&=\frac{-1}{4\pi}\int_{|\omega|=1}v(x+R\omega)d\omega -\frac1{4\pi
R}\int_{\partial B_R(x)}\frac{\partial}{\partial \nu}v(y)dS(y)\\
&-\int_{B_R}\frac1{4\pi|x-y|}\text{div}(\chi_\Omega(y)u(y))dy\\
&=\frac{-1}{4\pi}\int_{|\omega|=1}v(x+R\omega)d\omega-\int_{B_R}\frac1{4\pi|x-y|}\text{div}(\chi_\Omega(y)u(y))dy.
\endaligned$$
Notice that from $\nabla v\in L^2(\Bbb R^3)$, we have
$$
\lim_{R\to\infty}\nabla_x\int_{|\omega|=1}v(x+R\omega)d\omega=0.
$$
So
$$\aligned
\nabla v(x) &=\frac{-1}{4\pi}\nabla\int_{\Bbb R^3}\frac1{|x-y|}\text{div}(\chi_\Omega(y)u(y))dy\\
&=\frac1{4\pi}\nabla\int_\Omega u(y)\cdot\nabla\frac1{|x-y|}dy.
\endaligned$$

We have (c.f.[M]Th.2.6.2)
$$
\nabla v(x)=\Cal{L}( u)=Au(x)\chi_\Omega(x)+\frac1{4\pi}\int_\Omega\nabla_x(u(y)\cdot\nabla_x\frac1{|x-y|})dy,
$$
for $u\in C^\alpha_{loc}(\Omega)$, where $A$ is a constant matrix. It is easy to check that $\Cal{L}-A$ is a
 Calder\'{o}n-Zygmund operator ([MC] Chapter 7). Then it is bounded from
$$\aligned
&\Cal{H}^1\text{(Hardy)}\quad\text{to}\quad L^1(\Bbb R^3),\\
&L^\infty(\Omega) \quad\text{to}\quad \text{BMO},\\
&L^p(\Omega)\quad\text{to}\quad L^p(\Bbb R^3), \quad\forall p\in(1,\infty).
\endaligned$$
Thus the same is true for $\Cal{L}$.\qed

\bigskip

{\bf 3. Expressing Landau-Lifshitz-Maxwell equations in spherical coordinates}

Denote
$$u(x)=(\cos\xi(x)\cos\theta(x),\cos\xi(x)\sin\theta(x),\sin\xi(x)),\quad (-\frac\pi2\leqq \xi \leqq \frac\pi2,\,\, \theta \in S^1=\Bbb
R/2\pi\Bbb Z)
$$
and
$$
g(x)=(\cos\theta_g(x),\sin\theta_g(x),0).
$$

For convenience, we also use the notations
$$u(\theta,\xi)=(\cos\xi(x)\cos\theta(x),\cos\xi(x)\sin\theta(x),\sin\xi(x))$$
and
$$
u_\theta(\theta,\xi)=(-\cos\xi\sin\theta,\cos\xi\cos\theta,0),\quad u_\xi(\theta,\xi)=(-\sin\xi\cos\theta,-\sin\xi\sin\theta,\cos\xi)
$$
and etc.

By these notations, the energy functional $E_\lambda$  is rewritten as
$$
E_\lambda(\theta,\xi) =\int_\Omega
(\frac{1}2|\nabla\xi|^2+\frac{\cos^2\xi}2|\nabla
\theta|^2-\frac12 u(\theta,\xi)\cdot \nabla v+\lambda
\sin^2\xi) dx.\tag3.1
$$

The Euler-Lagrange
equation of (3.1) can be written as
$$\left\{\aligned
&\Delta\xi -
(\lambda-\frac{|\nabla\theta|^2}{2})\sin2\xi+u_\xi(\theta,\xi)\cdot\nabla v=0 \quad\text{in}\quad \Omega\\
&\xi=0 \quad\text{on}\quad \partial\Omega
\endaligned\right.\tag3.2
$$
and
$$\left\{\aligned
&\text{div}(\cos^2\xi
\nabla\theta)+u_\theta(\theta,\xi)\cdot\nabla v=0
\quad\text{in}\quad \Omega\\
&\theta=\theta_g \quad\text{on}\quad
\partial\Omega\\
&\Delta v=(-1)\text{div}\{u(\theta,\xi)\chi_\Omega\}\quad\text{on}\quad\Bbb R^3.
\endaligned\right.\tag3.3
$$

\bigskip

{\bf 4. Limit case}

Let $\lambda\to\infty$ in (3.1) and consider the limit functional
$$
E_\infty(\theta)=\int_\Omega\frac12|\nabla\theta|^2-\frac12u(\theta,0)\cdot\nabla v dx.
$$

Its critical points are the maps to $S^1$ which satisfy
$$\left\{\aligned
&\Delta\theta+u_\theta(\theta,0)\cdot\nabla v=0 \quad\text{in}\quad
\Omega\\
&\theta=\theta_g \quad\text{on}\quad
\partial\Omega\\
&\Delta v=(-1)\text{div}\{u(\theta,0)\chi_\Omega\}\quad\text{on}\quad\Bbb R^3.
\endaligned\right.\tag4.1
$$

\proclaim{Lemma 4.1} Assume that $\partial\Omega$ is uniformly
$C^2$. There is a constant $d_0>0$ such that if $\Omega$ is not
simply connected, $|\Omega|\leq d_0$,
 then in the homotopy
class $[\theta_g]$ of continuous mappings from $\overline\Omega$
into $S^1$ there exists a unique
 solution $\theta_*\in C^{1+\alpha}(\overline\Omega,S^1)$ and
 $\nabla v_*\in L^p(\Bbb R^3)\cap BMO(\Bbb R^3)$ ($1<\forall p<\infty$) to (4.1).

\endproclaim

$Proof.$ Step 1. Existence.

 For any given $\nabla v\in L^p(\Bbb R^3)\cap L^2(\Bbb R^3) (1<p<\infty)$, there is a minimizing sequence of $E_\infty$ which converges to a minimizer $\theta$ of $E_\infty$. Let $\Cal{\Theta}$ denote the map from $\nabla v\in L^p(\Bbb R^3)\cap L^2(\Bbb R^3)$ to $\theta$ by
 $$\left\{\aligned
&\Delta\theta+u_\theta(\theta,0)\cdot\nabla v=0 \quad\text{in}\quad
\Omega\\
&\theta=\theta_g \quad\text{on}\quad
\partial\Omega.
\endaligned\right.
$$
Note that $\Delta\theta\in L^p(\Omega)$ ($1< p<\infty$).
From [L] Theorem 3.1.1, $\theta\in W^{2,p}(\Omega)$ ($1<
p<\infty$). By the Sobolev embedding theorem ([Ad]), for any
$\alpha\in(0,1)$, $\theta\in C^{1+\alpha}(\overline{\Omega})$.

For any $\nabla v_1$ and $\nabla v_2\in L^p(\Bbb R^3)\cap L^2(\Bbb R^3)$, let $\theta_i=\Cal{\Theta}(\nabla v_i)$ ($i=1,2$). Then
$$\aligned
&\int_\Omega|\nabla(\theta_1-\theta_2)|^2dx\\
&=\int_\Omega(\theta_1-\theta_2)\{(u_\theta(\theta_1,0)-u_\theta(\theta_2,0))\cdot\nabla v_1+u_\theta(\theta_2,0)\cdot(\nabla v_1-\nabla v_2)\}dx\\
&\leqq \|\nabla v_1\|_{L^{3/2}(\Omega)}\|\theta_1-\theta_2\|_{L^6(\Omega)}^2+\|\nabla(v_1-v_2)\|_{L^2(\Omega)}\|\theta_1-\theta_2\|_{L^2(\Omega)}.
\endaligned\tag4.2
$$

From the Sobolev inequality and the Poinc\'{a}re inequality,
$$\aligned
&\|\theta_1-\theta_2\|^2_{L^6(\Omega)} \leq
C_1\|\nabla(\theta_1-\theta_2)\|_{L^2(\Omega)}^2,\\
&\|\theta_1-\theta_2\|_{L^2(\Omega)} \leq
C|\Omega|^{1/3}\|\nabla(\theta_1-\theta_2)\|_{L^2(\Omega)}
\endaligned\tag4.3$$
where the constants $C_1$, $C$ are independent of $\Omega$. See [GT], we may take
$$C_1=(\frac1{3\pi})^{1/2}(\frac3{\Gamma(5/2)})^{1/3}.$$

On the other hand, we have
$$\aligned
&\int_{\Bbb R^3}|\nabla v|^2dx=(-1)\int_\Omega u(\theta,\xi)\cdot\nabla vdx\\
&\leq |\Omega|^{1/2}(\int_\Omega|\nabla v|^2dx)^{1/2}.
\endaligned$$
So we have
$$
\|\nabla v\|_{L^2(\Bbb R^3)}\leq |\Omega|^{1/2}\tag4.4
$$
and
$$
\|\nabla v\|_{L^{3/2}(\Omega)}\leq |\Omega|^{1/6}\|\nabla
v\|_{L^2(\Omega)}\leq |\Omega|^{2/3}. \tag4.5$$
From (4.2),(4.3) and (4.5), there is $d_0>0$ such that if $|\Omega|\leq d_0$ then
$$
\int_\Omega|\nabla(\theta_1-\theta_2)|^2dx\leq C\|\nabla(v_1-v_2)\|_{L^2(\Omega)}^2
\tag4.6$$
which implies that $\Cal{\Theta}$ is a continuous map from $L^2(\Bbb R^3)\cap L^p(\Bbb R^3)$ to $C^{1+\alpha}(\overline{\Omega},S^1)$.

Let  $[\theta_g]$ denote the homotopy class of $\theta_g$ in
$C^\alpha(\overline{\Omega},S^1)$. Note that $[\theta_g]$ is
convex and closed. From Lemma 2.1, $\nabla v\in L^p(\Bbb R^3)\cap
BMO(\Bbb R^3)$ for $p\in(1,\infty)$, then the map
$\Cal{\Theta}\Cal{L}u(\theta,0)$ is continuous from $\theta\in
[\theta_g]$ to
$[\theta_g]$, and
$\Cal{\Theta}\Cal{L}u([\theta_g],0)$ is pre-compact in
$[\theta_g]\subset C^\alpha(\overline{\Omega},S^1)$. Then by the Schauder fixed
point theorem, we proved the existence of solutions to (4.1).

Step 2. Uniqueness.

Because
$$
\int_{\Bbb R^3} |\nabla(v_1-v_2)|^2dx =(-1)\int_\Omega (u(\theta_1,0)-u(\theta_2,0))\cdot\nabla(v_1-v_2)dx,
$$
we have
$$\aligned
&(\int_{\Bbb R^3} |\nabla(v_1-v_2)|^2dx)^{1/2}\leq (\int_\Omega|u(\theta_1,0)-u(\theta_2,0)|^2dx)^{1/2}\\
&\leq |\Omega|^{1/3}(\int_\Omega|u(\theta_1,0)-u(\theta_2,0)|^6dx)^{1/6}\leq 2C_1|\Omega|^{1/3}\|\nabla(\theta_1-\theta_2)\|_{L^2(\Omega)}.
\endaligned$$
From (4.6), we have
$$
(\int_\Omega|\nabla(\theta_1-\theta_2)|^2dx)^{1/2}\leq CC_1|\Omega|^{1/3}(\int_\Omega|\nabla(\theta_1-\theta_2)|^2dx)^{1/2}.
$$
So there is $d_0>0$ such that if $|\Omega|\leq d_0$, then $\theta_1=\theta_2$.\qed

\bigskip

\proclaim{Lemma 4.2} Assume that $\partial \Omega$ is uniformly
$C^4$. Then the solution obtained in Lemma 4.1 satisfies
$$
\nabla v_*\in C^{3+\alpha}(\overline{\Omega}),\quad \theta_*\in
C^{3+\alpha}(\overline{\Omega}). $$ Here the derivatives
on the boundary $\partial\Omega$ take the inner limit of the
derivatives in $\Omega$  respectively.
\endproclaim

$Proof.$ Since $u_*(x)=(\cos\theta_*(x),\sin\theta_*(x),0)\in
W^{1,p}({\Omega})$, by [Ad] Theorem 4.26, there is an extension
$U_*\in W^{1,p}( B_R)$ of $u_*$ for $B_R\supset\supset\Omega$, and
$$
\|U_*\|_{W^{1,p}( B_R)}\leq C\|u_*\|_{W^{1,p}(\Omega)}
$$
where $C$ is independent of $u_*$.

Consider the equation
$$
\Delta V=(-1)\text{div}(U_*\chi_{B_R})   \quad\text{in}\quad \Bbb R^3. \tag4.7$$
By
elliptic equation theory ([GT]), we have
$$
\|V\|_{W^{2,p}(\Bbb R^3)}\leq C\|U_*\|_{W^{1,p}( B_R)}\leq
C\|u_*\|_{W^{1,p}(\Omega)}.
$$
Since on $\Omega$, $(-1)\text{div }(U_*\chi_{B_R})=(-1)\text{div}(u_*\chi_\Omega)=\Delta v_*$, we have $\nabla v_*\in
W^{1,p}(\Omega)$, and from (4.1), $\theta_*\in W^{3,p}(\Omega)$
and so is $u_*$. Using [Ad] Theorem 4.26 and (4.7) again,
$$
\|U_*\|_{W^{3,p}(B_R)}\leq C \|u_*\|_{W^{3,p}(\Omega)},\quad V\in W^{4,p}(\Bbb R^3)
$$
so $\nabla v_*\in W^{4,p}(\Omega)$, and from the
Sobolev inequality, $\nabla v_*\in C^{3+\alpha}(\overline\Omega)$.
From (4.1), $\theta_*\in C^{3+\alpha}(\overline\Omega)$.\qed

\bigskip

{\bf 5. Proof of theorem 1}

 Let $\alpha_0 \in (0,1)$ and define
$$\aligned
M(\theta_*,\nabla v_*) = \{(\theta,\nabla v) &| \theta \in C^{1+\alpha_0}(\overline\Omega),\,\,
\theta|_{\partial\Omega}=\theta_g,\,\, \theta\in [\theta_g],\,\,\|\theta-\theta_*\|_{C^{1+\alpha_0}(\overline\Omega)}\leq 1,\\
&v\in V,\,\,  \nabla v\in C^{\alpha_0}(\overline\Omega),\,\,
\|\nabla v-\nabla v_*\|_{C^{\alpha_0}(\overline\Omega)}\leq 1  \}.
\endaligned$$

\bigskip

\proclaim{Lemma 5.1} Assume that $\partial\Omega$ is uniformly
$C^2$. For any given $(\theta,\nabla v) \in M(\theta_*,\nabla
v_*)$, there exists a continuous map $\Cal{\Xi}_\lambda$ such that
 $\xi_\lambda=\Cal{\Xi}_\lambda(\theta,\nabla v)$ is a solution to (3.2) which
 satisfies ,
$$\|\xi_\lambda\|_{W^{2,q}(\Omega)}\leq
C_1,\quad \forall q\in (1,\infty)   \tag5.1
$$
and
$$
\|\xi_\lambda\|_{W^{1,q}(\Omega)} \leq \frac {C_2}{\sqrt{\lambda}},
\tag5.2
$$
$$
 \|\xi_\lambda\|_{C^{\alpha}(\overline\Omega)} \leqq \frac {C_2}{\sqrt{\lambda}},
\tag5.3
$$
provided  $\lambda$ is large enough. Here the constants
$C_i=C_i(\|\theta\|_{C^{1+\alpha_0}(\overline\Omega)},\|\nabla v\|_{C^{\alpha_0}(\overline\Omega)})$ $(i=1,2)$ are
independent of $\lambda$.
\endproclaim

$Proof.$ Let $\eta:= \xi + \frac{C}\lambda$, where $C$ is a constant to be determined in the proof.
The equation for $\eta$ is written as
$$
\left\{\aligned
&(-1)\Delta\eta=-(\lambda-\frac{|\nabla\theta|^2}2)\sin2(\eta-\frac
C\lambda)+u_\xi(\theta,\eta-\frac{C}\lambda)\cdot\nabla v \quad\text{in}\quad
\Omega\\
&\eta=\frac{C}\lambda \quad\text{on}\quad
\partial\Omega.
\endaligned\right.\tag5.4
$$

Let
$$
F(\eta)=-(\lambda-\frac{|\nabla\theta|^2}2)\sin2(\eta-\frac
C\lambda)+u_\xi(\theta,\eta-\frac{C}\lambda)\cdot\nabla v.
$$
It is easy to check that there exists a constant
$C=C(\|\theta\|_{C^{1+\alpha_0}(\overline\Omega)},\|\nabla v\|_{C^{\alpha_0}(\overline\Omega)})$ and $\lambda_0
(>0)$ such that
$$F(0)\geq 0, \quad F(\frac{2C}\lambda) \leq 0$$
for $\lambda \geq \lambda_0$. From [A], there exists a
non-negative solution $\eta_\lambda$ to (5.4) which satisfies
$$ 0 \leq \eta_\lambda \leq \frac{2C}\lambda,$$
provided  $\lambda \geq \lambda_0$. Let
$\xi_\lambda=\eta_\lambda -\frac C\lambda$. Thus
$\xi_\lambda$ is a solution of (3.2) and
$$ -\frac C\lambda \leq \xi_\lambda \leq \frac C\lambda,\tag
5.5
$$
provided  $\lambda \geq \lambda_0$.

Rewrite (3.2) as
$$
-\Delta\xi + 2\lambda\xi
=\lambda(2\xi-\sin2\xi)+\frac{|\nabla\theta|^2}2\sin2\xi+u_\xi(\theta,\xi)\cdot\nabla v,
\tag 5.6$$
and use [L] Theorem 3.1.3 to obtain $\forall q\in(1,\infty)$
$$\aligned
&\|\xi\|_{W^{1,q}(\Omega)}\leq \frac C{\sqrt{\lambda}}
(\lambda\|2\xi-\sin2\xi\|_{L^q(\Omega)}+\|\nabla\theta\|_{L^q(\Omega)}^2+\|\nabla
v\|_{L^q(\Omega)}),\\
&\|\xi\|_{W^{2,q}(\Omega)}\leq
C(\lambda\|2\xi-\sin2\xi\|_{L^q(\Omega)}+\|\nabla\theta\|_{L^q(\Omega)}^2+\|\nabla
v\|_{L^q(\Omega)}).
\endaligned$$

Note that for any $\delta > 0$, there exists $\lambda(\delta) >
0$ such that
$$\|2\xi_\lambda-\sin2\xi_\lambda\|_{L^q(\Omega)} \leq \frac{\delta}{\lambda} \|\xi_\lambda\|_{L^q(\Omega)},
$$
provided  $\lambda \geq \lambda(\delta)$. So we have (5.1)-(5.2).

Using (5.2) for $q>3$ and the Sobolev embedding theorem, we get
(5.3) with $\alpha\leq 1-\frac 3q$. \qed

\proclaim{Lemma 5.2} Assume that $\partial\Omega$ is uniformly $C^4$. There exists $d_0>0$ such that if $|\Omega|\leq d_0$, then for
$\xi_\lambda=\Cal{\Xi}_\lambda(\theta,\nabla v)$ obtained in Lemma 5.1, there exists a continuous map $\bar\Cal{\Theta}$ such that
 $(\bar\theta,\nabla\bar v)=\bar\Cal{\Theta}(\xi_\lambda)$ is a solution to (3.3). Moreover, there is $\lambda_0>0$ such that
$$
\|\bar\theta-\theta_*\|_{C^{2+\alpha_0}(\overline\Omega)}
,\quad \|\nabla\bar v-\nabla v_*\|_{C^{1+\alpha_0}(\overline\Omega)}
$$
are bounded provided that $\lambda\geq\lambda_0$ and
$$
\|\bar\theta-\theta_*\|_{C^{1+\alpha_0}(\overline\Omega)} \to 0,\quad \|\nabla\bar v-\nabla v_*\|_{C^{\alpha_0}(\overline\Omega)}\to0,
 \quad\text{as}\quad \lambda \to \infty
\tag5.7
$$
uniformly for $(\theta,\nabla v) \in M(\theta_*,\nabla v_*)$.
\endproclaim

$Proof.$ The proof of existence part is similar to  Lemma 4.1. Moreover, as in the proof of Lemma 4.2, if $(\bar\theta,\nabla\bar v)\in C^{1+\alpha}(\overline\Omega)\times L^p(\Bbb R^3)$ $(1<p<\infty)$ is a solution to (3.3), then $u(\bar\theta,\xi_\lambda)\in W^{1,p}(\Omega)$, and for $B_R\supset\supset\Omega$, there is an extension $U\in W^{1,p}(B_R)$ of $u$. From the Maxwell equation
$$
\Delta V=(-1)\text{div}(U\chi_{B_R})\quad\text{in}\quad \Bbb R^3,
$$
we have $V\in W^{2,p}(\Bbb R^3)$. So $\nabla\bar v\in W^{1,p}(\Omega)$, and from (3.3), $\bar\theta\in W^{3,p}(\Omega)$. By the Sobolev embedding theorem, $\bar\theta\in C^{2+\alpha_0}(\overline\Omega)$. Noting that $u(\bar\theta,\xi_\lambda)\in W^{2,p}(\Omega)$, by using the Maxwell equation again, we have $V\in W^{3,p}(\Bbb R^3)$ and $\nabla\bar v\in W^{2,p}(\Omega)$. Then $\nabla\bar v\in C^{1+\alpha_0}(\overline\Omega)$.

So we only need to prove (5.7). From (3.3) and (4.1), we obtain the equations for $\bar\theta-\theta_*$:
$$\left\{\aligned
&\text{div}(\cos^2\xi_\lambda \nabla(\bar\theta-\theta_*))\\
&=-\nabla\cos^2\xi_\lambda\cdot
\nabla\theta_*
+(1-\cos^2\xi_\lambda)\Delta\theta_*+(u_\theta(\bar\theta,0)-u_\theta(\bar\theta,\xi_\lambda))\cdot\nabla\bar v\\
&+(u_\theta(\theta_*,0)-u_\theta(\bar\theta,0))\cdot\nabla\bar v-u_\theta(\theta_*,0)\cdot(\nabla\bar v-\nabla v_*) \quad\text{in}\quad \Omega\\
&\bar\theta-\theta_*=0 \quad\text{on}\quad
\partial\Omega
\endaligned\right.\tag5.8
$$
and the equation for $\bar v-v_*$:
$$
\int_{\Bbb R^3}|\nabla(\bar
v-v_*)|^2dx=(-1)\int_\Omega(u(\bar\theta,\xi_\lambda)-u(\theta_*,0))\cdot\nabla(\bar
v-v_*)dx. \tag5.9$$

Multiplying (5.8) by $\bar\theta-\theta_*$, by standard elliptic
equation theory as well as (5.9) and Lemma 5.1, as in the proof
of Lemma 4.1, we can prove that there exists $d_0>0$ such that if $|\Omega|\leq d_0$ then
$$
\|\nabla(\bar\theta-\theta_*)\|_{L^2(\Omega)}\to0
\quad\text{as}\quad \lambda\to \infty, \tag5.10$$
uniformly for
$(\theta,\nabla v) \in M(\theta_*,\nabla v_*)$. Using (5.9)
again, we have
$$
\|\nabla\bar v-\nabla v_*\|_{L^2(\Bbb R^3)}\to0\quad\text{as}\quad
\lambda\to \infty, \tag5.11$$
uniformly for $(\theta,\nabla v) \in
M(\theta_*,\nabla v_*)$.

As in the Lemma 4.1-4.2, we can prove that for $\lambda$ large enough,
$$
\|\bar\theta-\theta_*\|_{C^{2+\alpha_0}(\overline\Omega)}
,\quad \|\nabla\bar v-\nabla v_*\|_{C^{1+\alpha_0}(\overline\Omega)}
$$
are bounded uniformly for $(\theta,\nabla v) \in M(\theta_*,\nabla
v_*)$. So we have (5.7).\qed

From Lemma 5.1-5.3, we obtain

\proclaim{Lemma 5.3} $\bar\Cal{\Theta}\Cal{\Xi}_\lambda(M(\theta_*,\nabla v_*))$ is pre-compact in
$M(\theta_*,\nabla v_*)$ and $\bar\Cal{\Theta}\Cal{\Xi}_\lambda$ is continuous provided  $\lambda$ is
large enough.
\endproclaim

\bigskip

\proclaim{Lemma 5.4} For any given $(\theta,\nabla v)\in
M(\theta_*,\nabla v_*)$, if $\theta\in
C^{3+\alpha_0}(\overline\Omega)$ and $\nabla v\in
C^{1+\alpha_0}(\overline\Omega)$, then the solution
$\xi_\lambda=\Cal{\Xi}_\lambda(\theta,\nabla v)$  of (3.2)
obtained in Lemma 5.1  satisfies
$$\|\xi_\lambda\|_{C^{2+\alpha_0}(\overline\Omega)}\leq
C_3,   \tag5.12
$$
and
$$
\lim_{\lambda\to\infty}\|\xi_\lambda\|_{C^{2}(\overline\Omega)}=0,\tag5.13
$$
where the constant
$C_3=C_3(\|\theta\|_{C^{2+\alpha_0}(\overline\Omega)},\|\nabla
v\|_{C^{1+\alpha_0}(\overline\Omega)})$ is independent of
$\lambda$.

\endproclaim

$Proof.$ Take the derivative  $\partial_j=\frac{\partial}{\partial
x_j}$ in the equation (3.2), and rewrite it as
$$\aligned
-\Delta
\partial_j\xi+2\lambda\partial_j\xi=
&2\lambda(1-\cos2\xi)\partial_j\xi
+\nabla\theta\cdot\nabla\partial_j\theta
\sin2\xi+|\nabla\theta|^2(\cos2\xi)\partial_j\xi\\
&+\partial_ju_\xi(\theta,\xi)\cdot\nabla
v+u_\xi(\theta,\xi)\cdot\nabla\partial_jv\\
&=:f\quad\text{in}\quad\Omega\\
&\xi=0\quad\text{on}\quad\partial\Omega.
\endaligned\tag5.14$$
Applying [L] Theorem 3.1.3 to (5.14), we have for any
$q\in(1,\infty)$, there is $\lambda_0(q)>0$ such that for
$\lambda\geq\lambda_0$,
$$
\|\partial_j\xi\|_{W^{2,q}(\Omega)}\leq C\|f\|_{L^q(\Omega)}.
$$
By the Sobolev embedding theorem, we get (5.12). (5.2) and (5.12)
imply (5.13).\qed

\bigskip

$ Proof\, of\, Theorem 1.$ From Lemma 5.3 and Schauder fixed
point theorem, $\bar\Cal{\Theta}\Cal{\Xi}_\lambda$ has a fixed
point $(\theta_\lambda,\nabla v_\lambda)$ in $M(\theta_*,\nabla
v_*)$ for large $\lambda$. Then we obtain a solution
$(\theta_\lambda,\xi_\lambda)$ and $v_\lambda$ to (3.2)-(3.3)
which has the properties stated in Theorem 1 by Lemma 5.1-5.3 and
Lemma 5.4. \qed

\bigskip

\indent
{\bf 6. Proof of theorem 2}

The time developing Landau-Lifshitz-Maxwell equations (1.3) can be written as
$$\left\{\aligned
\partial_t\theta &=\frac1{\cos^2\xi}\{\text{div}(\cos^2\xi\nabla\theta)+u_\theta(\theta,\xi)\cdot\nabla v\}\\
&-\frac{\gamma}{\cos\xi}\{\Delta\xi+(\frac{|\nabla\theta|^2}2-\lambda)\sin2\xi+u_\xi(\theta,\xi)\cdot\nabla v\}\quad\text{in}\quad\Omega\times\Bbb R^+\\
\partial_t\xi &=\Delta\xi+(\frac{|\nabla\theta|^2}2-\lambda)\sin2\xi+u_\xi(\theta,\xi)\cdot\nabla v\\
&+\frac{\gamma}{\cos\xi}\{\text{div}(\cos^2\xi\nabla\theta)+u_\theta(\theta,\xi)\cdot\nabla v\}\quad\text{in}\quad\Omega\times\Bbb R^+\\
&\theta=\theta_g,\quad\xi=0,\quad\text{on}\quad\partial\Omega\times\Bbb R^+\\
&\text{div}\{\nabla v+\chi_\Omega u(\theta,\xi)\}=0,\quad\text{in}\quad\Bbb R^3\times\Bbb R^+.
\endaligned\right.\tag6.1$$

For simplicity, we  denote the solution
$(\theta_\lambda,\xi_\lambda)$ obtained in section 4 by $(\theta,\xi)$ and $$
\nabla v=\Cal{L}( u(\theta,\xi)).
$$
The linearized operator $A_\lambda$ of the right terms of (6.1) is written as
$$ A_\lambda=\left(\matrix a_{11} & a_{12}\\ a_{21} & a_{22}
\endmatrix\right),
$$where
$$\aligned
a_{11}&=
\Delta+\frac{\nabla\cos^2\xi\cdot\nabla}{\cos^2\xi}-\frac{\gamma\sin2\xi}{\cos\xi}\nabla\theta\cdot\nabla+\frac1{\cos^2\xi}\{u_{\theta\theta}(\theta,\xi)\cdot\nabla v\\
&+u_\theta(\theta,\xi)\cdot\Cal{L}( u_\theta(\theta,\xi)\cdot )\}-\frac{\gamma}{\cos\xi}\{u_{\xi\theta}(\theta,\xi)\cdot\nabla v+u_\xi(\theta,\xi)\cdot\Cal{L}( u_\theta(\theta,\xi)\cdot )\}
\endaligned$$
$$\aligned
a_{12}&= \frac2{\cos^2\xi}\nabla\xi\cdot\nabla\theta-\frac{2\sin\xi}{\cos\xi}\nabla\theta\cdot\nabla +\frac1{\cos^2\xi}\{u_{\theta\xi}(\theta,\xi)\cdot\nabla v+u_\theta(\theta,\xi)\cdot\Cal{L}( u_\xi(\theta,\xi)\cdot )\}\\
&+\frac{2\sin\xi}{\cos^3\xi}u_\theta(\theta,\xi)\cdot\nabla v-(\frac{\gamma\sin\xi}{\cos^2\xi})\{\Delta\xi+(\frac{|\nabla\theta|^2}2-\lambda)\sin2\xi+u_\xi(\theta,\xi)\cdot\nabla v\}\\
&-\frac{\gamma}{\cos\xi}(\Delta+(|\nabla\theta|^2-2\lambda)\cos2\xi)-\frac{\gamma}{\cos\xi}\{u_{\xi\xi}(\theta,\xi)\cdot\nabla v+u_\xi(\theta,\xi)\cdot\Cal{L}( u_\xi(\theta,\xi)\cdot)\}
\endaligned
$$
$$\aligned
a_{21} &=
\sin2\xi
\nabla\theta\cdot\nabla+(\frac{\gamma}{\cos\xi})(\cos^2\xi\Delta+\nabla\cos^2\xi\cdot\nabla)+u_{\xi\theta}(\theta,\xi)\cdot\nabla v\\
&+u_\xi(\theta,\xi)\cdot\Cal{L}( u_\theta(\theta,\xi)\cdot)+\frac{\gamma}{\cos\xi}\{u_{\theta\theta}(\theta,\xi)\cdot\nabla v+u_\theta(\theta,\xi)\cdot\Cal{L}( u_\theta(\theta,\xi)\cdot)\}
\endaligned$$
$$\aligned
a_{22}&=
\Delta +(|\nabla\theta|^2-2\lambda)\cos2\xi+u_{\xi\xi}(\theta,\xi)\cdot\nabla v+u_\xi(\theta,\xi)\cdot\Cal{L}( u_\xi(\theta,\xi)\cdot)\\
&+(\frac{\gamma}{\cos\xi})\{-\sin2\xi\Delta\theta-2\cos2\xi\nabla\xi\cdot\nabla\theta-\sin2\xi\nabla\theta\cdot\nabla\}+(\frac{\gamma\sin\xi}{\cos^2\xi})\text{div}(\cos^2\xi\nabla\theta)\\
&+\gamma u_\theta(\theta,0)\cdot\Cal{L}( u_\xi(\theta,\xi)\cdot)
\endaligned
$$

Decompose the operator $A_\lambda$ into $\bar A_\lambda$ and
the perturbation $G$:
$$\align
G &=A_\lambda-\bar A_\lambda\\
&=\left(\matrix \frac{\nabla\cos^2\xi\cdot\nabla}{\cos^2\xi}-\frac{\gamma\sin2\xi\nabla\theta\cdot\nabla}{\cos\xi} &
\left\{\aligned
&-\frac{\sin2\xi\nabla\theta\cdot\nabla}{\cos^2\xi}-\frac{\gamma\Delta}{\cos\xi}\\
&+\frac{2\gamma\lambda}{\cos\xi}(\cos2\xi+\sin^2\xi)
\endaligned\right\}\\
G_{21} & -\frac{\gamma\sin2\xi\nabla\theta\cdot\nabla}{\cos\xi}
\endmatrix\right)\tag6.2
\endalign$$
where
$$
G_{21}=a_{21}-u_\xi(\theta,\xi)\cdot\Cal{L}( u_\theta(\theta,\xi)\cdot)
$$

We consider the spectrum of the operator $A_\lambda$ as the
perturbation one of the operator $\bar A_\lambda$. As in [Z3][Z11],  following Proposition 5.1-5.3 can be proved in the same way.

\bigskip

\proclaim{Proposition 6.1} Let $T=\beta I-\bar A_\lambda$. For $\delta>0$,
there exist $\beta>0$,
$\lambda_0=\lambda_0(\delta)$, and for $\lambda\geqq \lambda_0$ there
exists $\gamma_0(\lambda)>0$ such that for $\lambda \geqq \lambda_0$
and $\gamma \in [0,\gamma_0]$, we have
$$
\|G\Phi\|_H\leqq
\delta^{1/2}((\beta+1)\|\Phi\|_H+\|T\Phi\|_H)\quad\text{for}\quad
\Phi\in D(T).
$$
That is, $G$ is $T$-bounded with $T$-bound $b: b\leqq \delta^{1/2}$.
\endproclaim

\bigskip

For $j=1,2$, let
$$Y_j=\{(\matrix\theta\\ \xi\endmatrix)\in W^{j,2}(\Omega,\Bbb R^2):\quad (\matrix\theta\\ \xi\endmatrix)=(\matrix\theta_g\\ 0\endmatrix)\,\,\text{on}\,\,\partial\Omega\}
$$
where $Y_j$ are endowed with the norm and inner product of $W^{j,2}(\Omega,\Bbb R^2)$. Note that $Y_j$ are Hilbert spaces. As in [Z3][Z11], we have

\proclaim{Proposition 6.2} There exist $d_0>0$, $\lambda_0>0$, and for $\lambda\geq\lambda_0$ there exists $\gamma_0(\lambda)>0$ such that
$$
A_\lambda:\,\, D(A_\lambda)(\subset L^2(\Omega,\Bbb R^2))\to L^2(\Omega,\Bbb R^2)
$$
 is a sectorial operator with $D(A_\lambda)=Y_2$, provided that $|\Omega|\leq d_0$,
$\lambda\geqq\lambda_0$ and $\gamma\in[0,\gamma_0(\lambda)]$. Moreover the norm of $Y_2$ is equivalent to the graph norm of $A_\lambda$
$$
\|(\matrix\theta\\\xi\endmatrix)\|_{Y_2}\sim\|(\matrix\theta\\\xi\endmatrix)\|_{L^2(\Omega,\Bbb R^2)}+\|A_\lambda(\matrix\theta\\\xi\endmatrix)\|_{L^2(\Omega,\Bbb R^2)},
$$
and for any $\mu\in \rho(A_\lambda)$, $(\mu-A_\lambda)^{-1}$ is compact.
\endproclaim

\bigskip

Let $\bar\mu_1(\lambda), \bar\mu_2(\lambda), \dots,\bar\mu_k(\lambda),\dots$ and
$$(\bar\phi_k^\lambda,\bar\psi_k^\lambda)\in W^{1,2}_0(\Omega)\times W^{1,2}_0(\Omega), \quad
\|\bar\phi_k^\lambda\|^2_{L^2(\Omega)}+\|\bar\psi_k^\lambda\|^2_{L^2(\Omega)}=1, \quad k=1,2,\dots$$
denote the eigenvalues and eigenfunctions of $\bar A_\lambda$,
respectively. Assume
$$
Re\bar\mu_1(\lambda)\leq Re \bar\mu_2(\lambda)\leq \dots\leq Re\bar\mu_k(\lambda)\leq\dots.
$$

\proclaim{Propostion 6.3} There exist $d_0>0$, $\lambda_0>0$, $\bar\gamma>0$ and $C>0$
 such that $Re\bar\mu_k(\lambda)\geqq -C$
 provided that $|\Omega|\leq d_0$,
$\lambda\geqq\lambda_0$ and $\gamma\in[0,\bar\gamma]$.
Moreover if there is $\{\lambda_j\}_j$  such that
$$
\limsup_{\lambda_j\to\infty}Re\bar\mu_k(\lambda_j)<\infty,
$$
then
$$
\limsup_{\lambda_j\to\infty}\lambda_j\int_\Omega (\bar\psi_k^{\lambda_j})^2dx<\infty.
$$
\endproclaim

\bigskip

The eigenvalue problem for $\bar A_\lambda$ can be written as
$$\left\{\aligned
\Delta\phi
&+\frac2{\cos^2\xi}(\nabla\xi\cdot\nabla\theta)\psi-(\frac{\gamma}{\cos\xi})(|\nabla\theta|^2(\cos2\xi)\psi)-(\frac{\gamma\sin\xi}{\cos^2\xi})(\Delta\xi+\frac{|\nabla\theta|^2}2\sin2\xi)\psi\\
&+\frac1{\cos^2\xi}\{u_{\theta\theta}(\theta,\xi)\cdot\nabla v\phi+u_\theta(\theta,\xi)\cdot\Cal{L}( u_\theta(\theta,\xi)\phi )\}\\
&-\frac{\gamma}{\cos\xi}\{u_{\xi\theta}(\theta,\xi)\cdot\nabla v\phi+u_\xi(\theta,\xi)\cdot\Cal{L}( u_\theta(\theta,\xi)\phi )\}\\
&+\frac1{\cos^2\xi}\{u_{\theta\xi}(\theta,\xi)\cdot\nabla v\psi+u_\theta(\theta,\xi)\cdot\Cal{L}( u_\xi(\theta,\xi)\psi)\}+\frac{2\sin\xi}{\cos^3\xi}u_\theta(\theta,\xi)\cdot\nabla v\psi\\
&-(\frac{\gamma\sin\xi}{\cos^2\xi})\{u_\xi(\theta,\xi)\cdot\nabla v\psi\}-\frac{\gamma}{\cos\xi}\{u_{\xi\xi}(\theta,\xi)\cdot\nabla v\psi+u_\xi(\theta,\xi)\cdot\Cal{L}( u_\xi(\theta,\xi)\psi)\}\\
&=-\bar\mu\phi \quad\text{in}\quad \Omega\\
\phi &=0 \quad\text{on}\quad
\partial\Omega ,
\endaligned\right.\tag6.3$$
and
$$\left\{\aligned
&\Delta\psi+(|\nabla\theta|^2-2\lambda)(\cos2\xi)\psi+(\frac{\gamma}{\cos\xi})(-\sin2\xi\Delta\theta-2\cos2\xi\nabla\xi\cdot\nabla\theta)\psi\\
&+(\frac{\gamma\sin\xi}{\cos^2\xi})\text{div}(\cos^2\xi\nabla\theta)\psi+u_\xi(\theta,\xi)\cdot\Cal{L}( u_\theta(\theta,\xi)\phi)\\
&+(u_{\xi\xi}(\theta,\xi)\cdot\nabla v)\psi+u_\xi(\theta,\xi)\cdot\Cal{L}( u_\xi(\theta,\xi)\psi)+\gamma u_\theta(\theta,0)\cdot\Cal{L}( u_\xi(\theta,\xi)\psi)\\
&=-\bar\mu\psi \quad\text{in}\quad \Omega\\
\psi &=0 \quad\text{on}\quad
\partial\Omega.
\endaligned \right.\tag6.4$$

\proclaim{Lemma 6.4} Suppose there is $\{\lambda_j\}_j$  such that
$$
\limsup_{\lambda_j\to\infty}Re\bar\mu_k(\lambda_j)<\infty.
$$
Then there exist $\lambda_0>0$ and $\bar\gamma>0$
such that
$$\|\bar\psi_k^{\lambda_j}\|_{C^\alpha(\overline\Omega)} \leq \frac
C{\sqrt{\lambda_j}},\tag6.5
$$
$$
\|\bar\psi_k^{\lambda_j}\|_{C^{2+\alpha}(\overline\Omega)}\leq
C,\tag6.6
$$
and
$$
\|\bar\phi_k^{\lambda_j}\|_{C^{2+\alpha}(\overline\Omega)}\leq
C,\tag6.7
$$
provided  $\lambda_j \geqq \lambda_0$ and $\gamma \in
[0,\bar\gamma]$, where the constant $C$  only depends on $k$ and the
$C^{2+\alpha}(\overline\Omega)$ norm of
$(\xi,\theta)$.
\endproclaim

$Proof.$ Step 1. Since $\phi$, $\psi\in L^2(\Omega)$, from Lemma 2.1 we have
$$
\Cal{L}( u_\theta(\theta,\xi)\phi),\,\,\Cal{L}( u_\xi(\theta,\xi)\psi)\in L^2(\Bbb R^3).
$$
From (6.3)-(6.4) and [L]Theorem 3.1.3, we get $\phi$, $\psi\in W^{2,2}(\Omega)$. Since $\Omega$ is regular, by the Sobolev embedding theorem, $\phi$, $\psi\in C^\alpha(\overline\Omega)$ for some $\alpha\in (0,1)$. So, from Lemma 2.1, we have for $p\in(1,\infty)$
$$
\Cal{L}( u_\theta(\theta,\xi)\phi),\,\,\Cal{L}( u_\xi(\theta,\xi)\psi)\in L^p(\Bbb R^3).
$$
By using (6.3)-(6.4) and [L] Theorem 3.1.3 again, we obtain
$$
\phi,\,\, \psi\in W^{2,p}(\Omega),\quad \forall p\in (1,\infty)
\tag6.8$$
and there is constant $C$ such that
$$
\|\phi\|_{W^{2,p}(\Omega)},\quad \|\psi\|_{W^{2,p}(\Omega)}\leq C,\quad\text{uniformly for $\lambda_j\geq\lambda_0$.}
$$

By the Sobolev embedding theorem,
$$
\phi,\,\, \psi\in C^{1+\alpha}(\overline\Omega).
$$
As in the proof of Lemma 4.2, we have
$$
\Cal{L}( u_\theta(\theta,\xi)\phi),\,\,\Cal{L}( u_\xi(\theta,\xi)\psi)\in C^{1+\alpha}(\overline\Omega).
\tag6.9$$
Considering the equations for $\partial_{x_i}\phi$ and $\partial_{x_i}\psi$, similarly we have $\phi$, $\psi\in C^{2+\alpha}(\overline\Omega)$ and
$$
\|\phi\|_{C^{2+\alpha}(\overline\Omega)},\quad  \|\psi\|_{C^{2+\alpha}(\overline\Omega)}
$$
are bounded uniformly for $\lambda_j\geq\lambda_0$.

 Step 2. Applying  [L] Theorem 3.1.3 to (6.4), we obtain that for
 all $q\in (1,\infty)$,
 there exist $\bar\gamma>0$ and
$\lambda_0>0$ such that for $\gamma\in [0,\bar\gamma]$ and
$\lambda_j\geq \lambda_0$,
$$\|\psi\|_{W^{1,q}(\Omega)}\leq
\frac{C}{\sqrt{\lambda_j}},\tag6.10
$$
where the constant $C$ is independent of $\lambda_j$, $\gamma$. By
the Sobolev embedding theorem we get (6.5).\qed

\proclaim{Lemma 6.5} Assume $|\Omega|\leq d_0$.  There exist
$\lambda_0>0$ and $\bar\gamma>0$, such that for $\lambda\geqq
\lambda_0$, $\gamma\in[0,\bar\gamma]$,
$$
Re\bar\mu_1(\lambda) \geq \mu_0(>0),
$$
  where $\mu_0$ is independent of $\lambda$.
\endproclaim

$Proof.$ If not, there is $\{\lambda_j\}_j$  such that
$$
\limsup_{\lambda_j\to\infty}Re\bar\mu_1(\lambda_j)<\mu_0.
$$

From Lemma 6.4, (6.3) and (6.4) converge to the
eigenvalue problem
$$\left\{\aligned
&\Delta \phi+u_{\theta\theta}(\theta_*,0)\cdot(\nabla v_*)\phi+(u_\theta(\theta_*,0)-\gamma(0,0,1)\cdot h =-\mu\phi \quad\text{in}\quad \Omega\\
&h =\Cal{L}( u_\theta(\theta_*,0)\phi),\quad (0,0,1)\cdot h=0\quad\text{in}\quad\Omega\\
&\phi=0 \quad\text{on}\quad
\partial\Omega.
\endaligned\right.\tag6.11
$$
Integrating over $\Omega$ and using
$$
-\int_\Omega \phi u_\theta(\theta_*,0)\cdot hdx=\int_{\Bbb R^3}|h|^2dx
$$
 we have
$$\aligned
&Re\mu\int_\Omega\phi^2dx\\
 &=Re\int_\Omega|\nabla\phi|^2dx-\phi^2u_{\theta\theta}(\theta_*,0)\cdot\nabla v_*dx\\
&+\int_{\Bbb R^3}|h|^2dx\\
&\geq\int_\Omega|\nabla\phi|^2dx-(\int_\Omega|\nabla v_*|^{3/2}dx)^{2/3}(\int_\Omega|\phi|^6dx)^{1/3}+\int_{\Bbb R^3}|h|^2dx\\
&\geq c_0\int_\Omega\phi^2dx
\endaligned$$
provided that $|\Omega|\leq d_0$ and $d_0$ is small enough, where $c_0>0$ is a constant. Taking $\mu_0=c_0/2$,
we get a contradiction. So we proved this lemma.\qed

$ Proof\, of\, Theorem 2.$ For fixed $\lambda \geqq \lambda_0$, we
can choose  $\gamma \in (0,\bar\gamma)$  such that $\gamma\lambda$
is small enough. By a perturbation argument (c.f. [K]), we have
that for $|\Omega|\leq d_0$ and for fixed $\lambda \geqq
\lambda_0$, there exists $\gamma_0(\lambda)>0$ such that for
$\gamma \in [0,\gamma_0(\lambda)]$, the  eigenvalues
$\{\mu_k(\lambda)\}_k$ of operator $A_\lambda$ have same behavior
as $\bar A_\lambda$.  Then
$$Re(\mu_1(\lambda)) \geqq \frac {\mu_0}2.
$$
By using the result of the Proposition 6.2 and [L] (p.295, the remark behind the proof of Theorem 8.1.1), we have the local existence of solutions to (6.1) with initial data in a neighborhood of the steady state solutions obtained in the Theorem 1. Moreover  by using [L] (Theorem 9.1.2), the solutions which satisfy the conditions in the Theorem 2
are exponentially asymptotically stable. \qed

\bigskip

{\bf Acknowledgement} This work is supported by NSFC NO.10571157.

\Refs

[Ad] Adams R.A., Sobolev spaces, Academic Press(1975).

[A] Amann H., On the existence of positive solutions of nonlinear
elliptic boundary value problems, Indiana Univ. Math. J.,  {\bf 21}, 125-146(1971).

[ABV] Anzellotti G., Baldo S. and Visintin A., Asymptotic behavior of
the Landau-Lifshitz model of ferromagnetism, Appl. Math. Optim.,
{\bf 23}, 171-192(1991).

[BTW] Ball J.M., Taheri A., Winter M., Local minimizers in micromagnetics and related problems, Calc. Var., {\bf14}, 1-27(2002).

[C] Campanato S., Generation of analytic semigroups  by elliptic
operators of second order in the H\"{o}lder spaces, Ann. Scuola
Norm. Sup. Pisa Cl Sci. (4)8, 495-512(1981).

[CSU] Chang N-H, Shatah J., Uhlenbeck K., Schr\"{o}dinger maps, Communications on Pure and Applied Math., Vol.LIII, 590-602(2000).

[FR] Fardoun A., Ratto A., Harmonic maps with potential, Calc. Var., {\bf5}, 183-197(1997).

[GS] Guo B., Su F., The global solution for Landau-Lifshitz Maxwell equations, J. Partial Diff. Eqs., {\bf14}, 133-148(2001).

[GSh] Gustafson S., Shatah J., The stability of localized solutions of Landau-Lifshitz equations, Communications on Pure and Applied Math., Vol.LV, 1-24(2002).

[GT] Gilbarg D. and Trudinger N., Elliptic partial differential
equations of second order, Springer, New York(1983).

[HL] Hang F-B., Lin F-H., Static theory for planar ferromagnets and antiferromagnets, Acta Math. Sin. (Engl. Ser.),  {\bf17}, no. 4, 541--580(2001).

[HS] Hubert A., Sch\"{a}fer R., Magnetic domains, Springer(1998).

[JK]James R.D., Kinderlehrer D., Frustration in ferromagnetic materials, Continuum Mech. Thermodyn., {\bf2}, 215-239(1990).

[JMZ] Jimbo S., Morita Y. and Zhai J., Ginzburg-Landau equation and
stable solutions in a nontrivial domain, Commun. In Partial
Differential Equations, {\bf20}(11\&12), 2093-2112 (1995).

[JZ] Jimbo S. and Zhai J., Ginzburg-Landau equation with magnetic
effect: non-simple-connected domains, J. Math. Soc. of Japan, {\bf50}(4), 663-684(1998).

[JZ2] Jimbo S. and Zhai J., Domain perturbation method and local minimizers to Ginzburg-Landau functional with magnetic effect, Abstract and Applied Analysis, {\bf5:2}, 101-112(2000).

[JZ3] Jimbo S. and Zhai J., Instability in a geometric parabolic equation on convex domain, J. Diff. Equ., {\bf188}(2), 447-460(2003).

[K] Kato T., Perturbation theory for linear operators, Springer
Verlag, New York(1966).

[KS] Kohn R.V., Sternberg P., Local minimizers and singular perturbations, Proc. Roy. Soc. Edin.A., Vol.111, 69-84(1989).

[L] Lunardi A., Analytic semigroups and optimal regularity in
parabolic problems, Birkh\"{a}user(1995).

[LL] Landau L.D., Lifshitz E.M., On the theory of the dispersion of
magnetic permeability in ferromagnetic bodies, Phys. Z. Sowjetunion 8,
(1935), Reproduced in Collected Papers of L.D. Landau, Pergamon, New
York, 101-114(1965).

[LU] Ladyzhenskaya O.A., Ural'tseva N.,Linear and quasilinear elliptic equations,Academic Press(1968).

[M] Morrey C.B., Calculus of variations, Springer(1966).

[MC] Meyer Y., Coifman R., Wavelets, Cambridge(1997).

[RS] Rivi\`{e}re T., Serfaty S., Compactness, kinetic formulation, and entropies for a problem related to micromagnetics, Comm. Partial Differential Equations, {\bf28}, 249-269(2003).

[SC] Slodi\v{c}ka M., Cimr\'{a}k I., Nunerical study of nonlinear ferromagnetic materials, Appl. Numer. Math., {\bf46}(1), 95-111(2003).

[V] Visintin A., On Landau-Lifshitz' equations for ferromagnetism,
Japan J. Appl. Math., {\bf2}, 69-84(1985).

[Z1] Zhai J., Heat flow with tangent penalization converges to mean
curvature motion, Proc. Roy. Soc. Edinburgh Sect. A, {\bf128}(A), 875-894(1998).

[Z2] Zhai J., Heat flow with tangent penalization, Nonlinear Analysis, TMA, {\bf28}, 1333-1346(1997).

[Z3] Zhai J., Non-constant stable solutions to Landau-Lifshitz equation, Calc. Var. \& PDE, {\bf7}, 159-171(1998).

[Z4] Zhai J., Existence and behavior of solutions to the Landau-Lifshitz equation, SIAM J. Math. Anal., Vol.30, No.4, 833-847(1999).

[Z5] Zhai J., Dynamics of domain walls in ferromagnets and weak ferromagnets, Physics Letters A, {\bf234}, 488-492(1997).

[Z6] Zhai J., Theoretical velocity of domain wall motion in ferromagnets, Physics Letters A, {\bf 242}, 266-270(1998).

[Z7] Zhai J., Velocity of domain wall in ferromagnets with demagnetizing field, Physics Letters A, {\bf279}, 395-399(2001).

[Z8] Zhai J., Fang J.P., Full velocity of micromagnetic domain walls, Physics Letters A, {\bf318}, 137-140(2003).

[Z9] Zhai J., Fang J.P., Li L.J., Wave Map with Potential  and  Hypersurface Flow, Dynamical Systems And Differential Equations, Supplement Volume, 940-946(2005).

[Z10] Zhai J., Li L.J., New results on Landau-Lifshitz ferromagnets model,  Nonlinear Analysis, TMA, 63,5-7, e11-e21(2005).

[Z11] Zhai J., Stable solutions to Landau-Lifshitz-Maxwell equations, Indiana Univ. Math. J. 54 No.6, 1635-1660(2005).

\endRefs

\end